\magnification = 1200
\input amssym.def
\input amssym.tex
\def \qed {\hfill $\square$}
\def \R {\Bbb R}

\def \d {\partial}
\def \ms {\medskip}
\def \msi {\medskip\noindent}
\def \ssi {\smallskip\noindent}

\def \bsi {\bigskip\noindent}
\def \pari {\par\noindent}

\def \sm {\setminus }

\def \wt {\widetilde }

\overfullrule=0pt
\def \dist{\mathop{\rm dist}\nolimits}

\def \sup{\mathop{\rm sup}\nolimits}
\def \inf{\mathop{\rm inf}\nolimits}

\centerline{APPROXIMATION OF A REIFENBERG-FLAT SET}
\par
\centerline{BY A SMOOTH SURFACE}

\bigskip
\centerline{Guy DAVID}

\bigskip

\vskip 1cm
\noindent
{\bf R\'{e}sum\'{e}.}
On montre que si l'ensemble $E \i \R^n$ est bien approch\'{e}, \`{a}
l'\'{e}chelle $r_0$, par des plans de dimension $d$,
il existe une surface lisse $\Sigma_0$ de dimension $d$, 
qui est proche de $E$ \`{a} l'\'{e}chelle $r_0$.
Quand $E$ est Reifenberg-plat,
ceci permet d'appliquer un r\'{e}sultat de G. David et T. Toro 
[Memoirs of the AMS 215 (2012), 1012], 
et de montrer que $E$ est l'image de $\Sigma_0$  par un 
hom\'{e}omorphisme bi-H\"old\'{e}rien de $\R^n$. Si de plus
$d=n-1$ et $E$ est compact et connexe, alors $\Sigma_0$ est  
orientable, et $\R^n \sm E$ a exactement deux composantes connexes
que la construction ci-dessus permet d'approximer par 
l'int\'{e}rieur par des domaines lisses.

\bigskip \noindent  
{\bf Abstract.}
We show that if the set $E \i \R^n$ is well approximated 
at the scale $r_0$ by planes of dimension $d$, we can find 
a smooth surface $\Sigma_0$ of dimension $d$ which is close to $E$ 
at the scale $r_0$. When $E$ is a Reifenberg flat set, this allows us 
to apply a result of G. David and T. Toro 
[Memoirs of the AMS 215 (2012), 1012], 
and get a bi-H\"older homeomorphism of $\R^n$ that sends $\Sigma_0$ to $E$.
If in addition $d=n-1$ and $E$ is compact and connected, 
then $\Sigma_0$ is orientable, and $\R^n \sm E$ has exactly 
two connected components, which we can approximate from the 
inside by smooth domains.


\bigskip\noindent
{\bf 1. Introduction.}
\ms
The initial purpose of this paper is to prove that if $E \i \R^n$ is a 
set which, in every ball of radius $r_0$, is sufficient 
close to a $d$-dimensional plane, we can
find a smooth surface $\Sigma_0$ of the same dimension, which
is close to $E$ at the scale $r_0$.

When $E$ is a Reifenberg flat set of dimension $d$ 
(which means that the property above holds at all scales $r \leq r_0$),
this allows us to apply one of the main results of [DT], 
and get a bi-H\"older homeomorphism of the ambient space $\R^n$
that sends $\Sigma_0$ to $E$. In the present context, the
difference between the result of [DT] that we want 
to use and the classical topological disk theorem of [R] 
is not large, it only resides in the fact that we authorize $E$
to look like a smooth surface at the scale $r_0$, rather than just a 
plane.

So the first main result of this paper can be seen a small preparation
lemma that slightly weakens the assumptions in [DT]. 

In the special case of a compact connected Reifenberg-flat
set $E$ of codimension $1$, we will deduce from this that
$\R^n \sm E$ has exactly two connected components $V_1$ and $V_2$, 
so that one does not need to mention this
separation property as an additional assumption,
and we will construct smooth domains that approximate $V_1$ and $V_2$
from the inside; see Theorem 1.18. 

\ms
Let us now state the results more precisely.
Fix integers $0 < d < n$, let $E$ be a (nonempty) closed set in $\R^n$, 
and define the (bilateral P. Jones) numbers $\gamma(x,r)$, 
$x\in E$ and $0 < r < +\infty$, by
$$
\gamma(x,r) = \inf\big\{ d_{x,r}(E,P) \, ; \, 
P \in {\cal P}(x) \big\},
\leqno (1.1)
$$
where ${\cal P}(x)$ denotes the set of $d$-dimensional
affine planes that contain $x$, and 
$$\eqalign{
d_{x,r}(E,P) &= {1 \over r}\,\sup\big\{ \dist(y,P) \, ; \,
y\in E \cap B(x,r) \big\}
\cr&\hskip 2 cm
+ {1 \over r}\, \sup\big\{ \dist(y,E) \, ; \, y\in P \cap B(x,r) \big\}
}\leqno (1.2)
$$
is a normalized local Hausdorff distance from $E$ to $P$.

We shall assume that there is a radius $r_0 > 0$ such that
$$
\gamma(x,r_0) \leq \varepsilon
\ \hbox{ for } x\in E,
\leqno (1.3)
$$
and prove that if $\varepsilon$ is small enough, depending only on
$n$ and $d$, there is a smooth $d$-dimensional surface
$\Sigma_0 \i \R^n$ with no boundary, such that
$$
\dist(x,\Sigma_0) \leq C_0 \varepsilon r_0
\ \hbox{ for } x\in E
\leqno (1.4)
$$
and 
$$
\dist(x,E) \leq C_0 \varepsilon r_0
\ \hbox{ for } x\in \Sigma_0,
\leqno (1.5)
$$
where $C_0$ depends only on $n$ and $d$. 
More precisely, we will show that there is a
$\Lambda > 0$, that depends only on $n$ and $d$, 
such that for each $y\in \Sigma_0$, 
$$
\hbox{$\Sigma_0$ coincides with a smooth 
$\Lambda \varepsilon$-Lipschitz graph in 
$B(y,10^{-2} r_0)$.} 
\leqno (1.6)
$$
That is, we can find a $d$-plane $P_y$ through $y$ and a 
$\Lambda \varepsilon$-Lipschitz mapping $F_y : P_y \to P_y^\perp$ such that, if
$$
{\cal G}(F_y) = \big\{ w + F_y(w) \, ; \, w\in P_y \big\}
\leqno (1.7)
$$
denotes the graph of $F_y$, then 
$$
\Sigma_0 \cap B(y,10^{-2} r_0) 
= {\cal G}(F_y) \cap B(y,10^{-2} r_0).
\leqno (1.8)
$$
In addition, the $F_y$ are smooth, and there exist constants
$\Lambda_k$, $k \geq 1$, such that
$$
||D^k F_y||_\infty \leq \Lambda_k \varepsilon r_0^{1-k}.
\leqno (1.9)
$$
(For $k=1$ we already knew this, with $\Lambda_1 = \Lambda$.)
Let us summarize all this officially.

\ms\proclaim Theorem 1.10.
There exist constants $\varepsilon_0 > 0$, $C_0 \geq 1$, and $\Lambda_k$, 
$k \geq 1$, that depend only on $n$ and $d$, such that 
if $E \i \R^n$ is a nonempty closed set such that
(1.3) holds for some $r_0 > 0$ and some $\varepsilon \in (0,\varepsilon_0)$, 
then we can find a smooth $d$-dimensional surface
$\Sigma_0$ with the properties (1.4)-(1.9).

\ms
Notice that we do not require $E$ to be flat at all scales
smaller than $r_0$ (as will be done for the next two results). 

Theorem 1.10 is designed so that we can apply Theorem 12.1
in [DT], and get the following statement, where we 
decided to work with $r_0=1$ for convenience.

\ms\proclaim Corollary 1.11.
There exist constants $\varepsilon_1 \leq \varepsilon_0$
and $C_1 > 1$, that depend only on $n$ and $d$, such that if 
$E \i \R^n$ is a nonempty closed set such that
$$
\gamma(x,r) \leq \varepsilon
\ \hbox{ for $x\in E$ and $0 <r \leq 1$},
\leqno (1.12)
$$
and if $\Sigma_0$ is the smooth surface provided by
Theorem 1.10, then there is a bijective mapping $g: \R^n \to \R^n$
such that
$$
g(x)=x
\ \hbox{ when } \dist(x,E) \geq 10^{-1},
\leqno (1.13)
$$
$$
|g(x)-x| \leq C_1 \varepsilon
\ \hbox{ for } x\in \R^n,
\leqno (1.14)
$$
$$
C_1^{-1} |x'-x|^{1+C_1\varepsilon} \leq |g(x)-g(x')| 
\leq C_1 |x'-x|^{1-C_1\varepsilon}
\leqno (1.15)
$$
for $x, x' \in \R^n$ such that $|x'-x| \leq 1$,
and
$$
g(\Sigma_0) = E.
\leqno (1.16)
$$

\ms
Let us check this. 
We want to apply Theorem 12.1 in [DT] 
to the set $E' = A E$, where we choose 
$A = 2\cdot 10^4$, and with the open set
$$
U = \big\{ x\in \R^n \, ; \, \dist(x,E') > 3 \big\}.
\leqno (1.17)
$$
The assumption (12.1) of [DT], 
relative to the smooth set $\Sigma'_0 = A \Sigma_0$,
is satisfied because of (1.6)-(1.9)
(although only with the constant $C\varepsilon$, which
does not matter). The proximity condition (12.3)
follows from (1.4) and (1.5), and so Theorem 12.1 in [DT] 
gives a mapping $g'$ that satisfies (12.4)-(12.8) in [DT]. 

Here we take $g(x) = A^{-1} g'(Ax)$, (1.13) holds
because (12.4) in [DT] 
says that $g'(x) = x$ when $\dist(x,U) \geq 13$,
while (1.14) and (1.15) easily follow from
(12.6) and (12.7) in [DT]. Finally, (12.8) in [DT] 
says that $E' \cap U = g'(\Sigma'_0 \cap U)$. 
But $E' = E' \cap U$ by (1.17),
and $\Sigma'_0 \i U$ by (1.5), so in fact
$E' = g'(\Sigma'_0)$ and $E = g(\Sigma_0)$, as needed.
So Corollary 1.11 follows
from Theorem 1.10 and Theorem 12.1 in [DT]. 
\qed 

\ms
We shall also apply Theorem 1.10 in the special case of 
connected sets of codimension~$1$, and obtain smooth approximations
from the inside for each of the two complementary components.

\ms\proclaim Theorem 1.18.
Suppose $d=n-1$.
There exist constants $\varepsilon_2 \leq \varepsilon_0$,
$C_2 > 1$, and $\Lambda'_k$, $k \geq 1$, that depend only on $n$, 
such that if $E \i \R^n$ is a nonempty compact connected set such that
for some choice of $r_0 > 0$ and $\varepsilon \in (0,\varepsilon_0)$,
$$
\gamma(x,r) \leq \varepsilon
\ \hbox{ for $x\in E$ and $0 < r \leq r_0$,}
\leqno (1.19)
$$
then
$$
\hbox{$\R^n \sm E$ has exactly two connected components,}
\leqno (1.20)
$$
which we shall denote by $V_1$ and $V_2$. In addition, for 
$0 < r < r_0$, there exist smooth connected domains
$W_{r,1} \i V_1$ and $W_{r,2} \i V_2$, with the following properties:
$$
\R^n \sm [W_{r,1} \cup W_{r,2}] \i \big\{ x\in \R^n \, ; \,
\dist(x,E) \leq 2 C_2 \varepsilon r \big\}
\leqno (1.21)
$$
and, for $j = 1, 2$,
$$
C_2 \varepsilon r \leq \dist(x,E) \leq 2 C_2 \varepsilon r
\ \hbox{ for } x\in \d W_{r,j} 
\leqno (1.22)
$$
and we have a local Lipschitz description of $\d W_{r,j}$ as in 
(1.6)-(1.9). That is, for each $y\in \d W_{r,j}$,
we can find a hyperplane $ P_{r,j,y}$ through $y$
and a Lipschitz mapping
$F_{r,j,y} : P_{r,j,y} \to P^\perp_{r,j,y}$ such that
$$
\hbox{$\d W_{r,j}$ coincides with the
graph of $F_{r,j,y}$ inside $B(y,10^{-3}r)$}
\leqno (1.23)
$$
and, for $k \geq 1$,
$$
||D^k F_{r,j,y}||_\infty \leq \Lambda'_k \varepsilon r^{1-k}.
\leqno (1.24)
$$
Finally, $W_{r,j} \cap B(y,10^{-3}r)$ lies on only one side of 
the graph of $F_{r,j,y}$.

\ms
Notice that (1.20) is part of the conclusion: as we shall see,
the Reifenberg-flatness of $E$ implies that it separates, and the 
fact that there are at most two complementary components follows from its
connectedness.

When $E$ is a little more regular than just Reifenberg-flat, we can 
expect to find other, more precise ways of saying that the $W_{r,j}$
converge to $V_j$, but hopefully the main ingredient will be 
(1.22)-(1.24) and the regularity of $E$. For instance, if $E$ is
also Ahlfors-regular of dimension $n-1$, so are the $\d W_{r,j}$, 
with uniform bounds, and the characteristic function of $W_{r,j}$
converges, locally in $BV(\R^n)$, to the characteristic function of $V_j$.

\ms
Let us give a short proof of (1.20) that relies on Corollary 1.11.
We shall see in Section~3 another proof that does not use the
parameterization from [DT], but only Theorem 1.10 and some small
amount of tracking at different scales.

Let $E$ be as in the statement. By scale invariance and without loss of 
generality, we may assume that $r_0 = 1$, and this way we will apply
Corollary 1.11 without rescaling.
Let $\Sigma_0$ and $g$ be as in Theorem 1.10 and then Corollary 1.11. 
We know that $\Sigma_0$ is a smooth hypersurface without boundary, 
and it is compact by (1.5). We shall soon check that 
$$
\Sigma_0 \hbox{ is connected;}
\leqno (1.25)
$$
let us conclude from here.
By [AH], page 440 (for nonsmooth manifolds), or rather [S] 
(a shorter and intuitive argument for the smooth case, using 
transversality), $\Sigma_0$ is orientable and $\R^n \sm \Sigma_0$ 
has exactly two connected components, which we denote by $U_1$ and $U_2$. 
Recall that $g$ is smooth, bijective, and equal to the
identity near infinity; then we get that
$\R^n \sm E = g(\R^n \sm \Sigma_0)$ has two components, 
namely the $g(U_j)$, $j=1,2$.

We still need to check that $\Sigma_0$ is connected. 
Let $a_1, a_2 \in \Sigma_0$ be given.
By (1.5), we can find $b_1, b_2 \in E$ such that 
$|b_i-a_i| \leq C \varepsilon$. Because $E$ is connected, 
we can find a chain of points $w_i \in E$, $0 \leq i \leq m$, with
$w_0 = b_1$, $w_m = b_2$, and $|w_i-w_{i-1}| \leq C \varepsilon$
for $1 \leq i \leq m$. 
(Otherwise, the set of points of $E$ that can be connected to $b_1$
by such a chain, which is open and closed in $E$, would contain 
$b_1$ but not $b_2$.) By (1.4), we can find $y_i \in \Sigma_0$
such that $|y_i-w_i| \leq C \varepsilon$. By definition
of $b_1$ and $b_2$, we can take $y_0 = a_1$ and $y_m = a_2$. 
Finally, by (1.6)-(1.8), applied to the points $y_i$, we get 
that for $1 \leq i \leq m$, there is a curve in $\Sigma_0$ 
that connects $y_{i-1}$ to $y_i$.
Thus $\Sigma_0$ is connected, and (1.20) follows.
\qed

\msi{\bf Remark 1.26.}
Instead of assuming that $E$ is connected, it is enough to assume that
it is ${r_0 \over 20}$-connected.
That is, that any two points $b_1$, $b_2$ in $E$ can be connected
by a chain $\{w_i \}$ in $E$ as above, with 
$|w_i-w_{i-1}| \leq {r_0 \over 20}$ for $1 \leq i \leq m$.
Indeed, our proof of (1.20) shows that $\Sigma_0$ is connected,
and then $E = g(\Sigma_0)$ is connected too. 
We can even prove that $E$ is connected without using Corollary 1.11;
see the proof of (3.1) below.

If we do not suppose that $E$ is connected (or  ${r_0 \over 
20}$-connected) we only get that
$\R^n \sm E$ has at least two components, but we still
have the approximation result for each component. This is not
hard: let $V$ be any component of $\R^n \sm E$, 
pick any $x_0 \in E$ such that $\dist(x_0,V) < r_0/100$,
and apply the result to the set $E'$ of 
points $x\in E$ that can be ${1\over 20}$-connected to
$x_0$. By Theorem 1.18, $E'$ separates $\R^n$ into two components,
and our approximation result holds for those. It is easy
to see that one of them contains $V$, and this gives
a good description of $V$ near $x_0$.

\msi{\bf Remark 1.27.}
We do not really need to assume that $E$ is bounded either;
Theorem 1.18 is still valid when $E$ is closed (instead of compact) 
and we keep the same other assumptions. The only place where we need
to change the proof is when we apply [S], to show that smooth compact
sets are orientable and have exactly two complementary components,
but the transversality argument still works for the surfaces that
we construct. It is probable that the results of this paper 
can also be extended to unbounded sets, even when the scale $r_0$ 
under which (1.19) is assumed to hold depends gently on on $x$. 
But we shall not try to do this.

\msi{\bf Remark 1.28.}
In [DT] the authors also consider situations 
where instead of requiring the bilateral control (1.12),
one merely assumes that the points of $E$ lie close to
$d$-planes $P(x,r)$, with some control on how fast they 
depend on $(x,r)$. Here we can try to construct $\Sigma_0$
with similar data; our construction will only give 
a smooth surface $\Sigma_0$, but with a boundary, and we
shall not try to see whether it is contained in a smooth 
surface without boundary.
See Remark 2.31. 

\ms
Theorem 1.10 will be proved in Section 2, by a 
simple construction where we start from a net of points 
and use a covering to fill the holes in a finite number of steps.

For Theorem 1.18 (proved in Section 3), we will first apply 
Theorem 1.10 at all intermediate scales $0<r \leq r_0$ 
to get central surfaces $\Sigma_{r}$, then move in the normal 
direction (in both directions, and until we are far enough from $E$) 
to get two surfaces $\Sigma_{r,1}$ and $\Sigma_{r,2}$, and show that these 
surfaces bound domains $W_{r,1}$ and $W_{r,2}$ that satisfy the
desired properties. The control of the connected components, if not 
surprising, will be the most unpleasant part of the argument. It
will be taken care of by a top-down compatibility argument.

\ms
The author wishes to thank T. D. Luu for discussions about 
Theorem 1.10 (we thought at some time that a result like this
would be needed for his PhD thesis, but it  turned out not to be the
case), T. Toro for suggesting that Theorem 1.10 be used to approximate
Reifenberg-flat domains from inside, and A. Lemenant and A. Chambolle, 
in particular for discussions about (1.20) and
orientability. He wishes to acknowledge support from the
Institut Universitaire de France, and the grant 
ANR-12-BS01-0014-01 "GEOMETRYA".

\medskip\noindent
{\bf 2. A proof of Theorem 1.10.}
\ms
In this section we prove Theorem 1.10. Since the
statement is invariant under dilations, it will be enough 
to prove the theorem when $r_0 = 1$. 

Let $E$ be as in the theorem, with $r_0 = 1$, set
$a = {1 \over 32}$, and choose a maximal set
$X \i E$ such that $|x-x'| \geq a$ for
$x, x' \in X$ such that $x'\neq x$. Thus
$$
\dist(x,X) \leq a
\ \hbox{ for } x\in E,
\leqno (2.1)
$$
because otherwise we could add $x$ to $X$. Next decompose
$X$ as a disjoint union
$$
X = \bigcup_{1 \leq j \leq N} X_j,
\leqno (2.2)
$$
where for each $j$, 
$$
|x-x'| \geq 16a = {1 \over 2} 
\ \hbox{ when $x$ and $x'$ lie in $X_j$ and } x \neq x'.
\leqno (2.3)
$$
Let us check that we can do this with an $N$ that depends only on $n$.
Define $X_j$ by induction, to be a maximal set
contained in $X \sm \cup_{i < j} X_i$ and with the
property (2.3); it is easy to see that we can 
stop as soon as $j$ is larger than the maximal number 
of points in a ball of radius $16a$ that lie at mutual
distances at least $a$. This number is in turn estimated
by saying that the balls of radius $a/2$ centered at these points
are all disjoint and contained in a ball of radius $17a$, which yields
$N \leq 34^n$ by comparing Lebesgue measures.

For each $x\in X$, use (1.3) to find an affine 
$d$-plane $P_x$ through $x$ such that
$$
d_{x,1}(E,P_x) \leq \varepsilon
\leqno (2.4)
$$
(see the definitions (1.1) and (1.2)). 
Then denote by $\pi_x : \R^n \to P_x$ the orthogonal 
projection on $P_x$, and by $\pi_x^\perp$ the orthogonal 
projection on the vector space $P_x^\perp$ of dimension
$n-d$ which is orthogonal to $P_x$.
Let us check that 
$$
d_{x,1/4}(P_x,P_y) \leq 8 \varepsilon
\ \hbox{ for $x, y \in X$ such that } |x-y| \leq {1 \over 2}.
\leqno (2.5)
$$
If $z\in P_x \cap B(x,1/4)$, we can use
(2.4) to find $w\in E$ such that $|w-z| \leq \varepsilon$;
then $w\in E \cap B(y,1)$ and (2.4) for $y$ says that
we can find $z'\in P_y$ such that $|z'-w| \leq \varepsilon$
and hence $|z'-z| \leq 2\varepsilon$.
By a similar argument, for each $z'\in P_y \cap B(x,1/4)$
we can find $z\in P_x$ such that $|z-z'| \leq 2\varepsilon$,
and (2.5) follows.

Thus in (2.5), $P_x$ and $P_y$ make a small angle;
this will be useful because we want small Lipschitz
graphs over $P_x$ to be small Lipschitz graphs over $P_y$ 
as well.

We shall now construct a nondecreasing sequence of sets
$S_j$, $0 \leq j \leq N$. We start with $S_0 = X$,
and our final set $S_N$ will be a good choice of
$\Sigma_0$, except for the fact that we shall not
immediately take care of the many derivatives in (1.9).
Notice that $X$ gives the general position
of $\Sigma_0$, so our problem will essentially consist in
completing $S_0$ into a smooth surface that contains it; 
we shall only use the set $E$ marginally,
to prove that $\Sigma_0$ has no boundary.

We shall construct the $S_j$ by induction,
with the property that for each $x\in X$,
there is an $A_j \varepsilon$-Lipschitz mapping
$F_{j,x} : P_x \to P_x^\perp$ such that
$$
S_{j} \cap B(x,4a) \i {\cal G}(F_{j,x}),
\leqno (2.6)
$$
where ${\cal G}(F_{j,x})$ denotes the graph of
$F_{j,x}$ over $P_x$, defined as in (1.7). 
The constants $A_j$ will be chosen larger and
larger, but since we can take $\varepsilon$ as 
small as we want, the $A_j \varepsilon$ will stay small.

In order to prove (2.6), we shall first check that
for each $x\in X$,
$$
|\pi_x^\perp(y)-\pi_x^\perp(z)| 
\leq A'_j \varepsilon |\pi_x(y)-\pi_x(z)|
\ \hbox{ for } y, z \in S_{j} \cap B(x,4a)
\leqno (2.7)
$$
for some constant $A'_j$ (that will depend 
on $A_{j-1}$ if $j \geq 1$).
As soon as we have (2.7), we observe that $\pi_x$
is injective on $S_{j} \cap B(x,4a)$, so
can define a function $F_{j,x}$ from
$H= \pi_x(S_{j} \cap B(x,4a))$ to $P_x^\perp$, 
by the relation $F_{j,x}(\pi_x(y)) = \pi_x^\perp(y)$
for $y\in S_{j} \cap B(x,4a)$. In addition, (2.7) says that
$F_{j,x}$ is $A'_j \varepsilon$-Lipschitz on $H$, and 
$S_{j} \cap B(x,4a)$ is its graph.

We then use the simpler Whitney extension theorem
(see [St], page 170) to extend $F_{j,x}$ 
into an $A_j \varepsilon$-Lipschitz function defined on $P_x$, 
and we get (2.6). Since we allow ourselves to take $A_j$ larger 
than $A'_j$, we don't need to use Kirszbraun's theorem 
(see 2.10.43 in [F]) 
and the construction of the extension is simpler. 
o
Let us now check that (2.7) holds for $j=0$. Recall
that we took $S_0 = X$.
Let $y, z \in S_{0} \cap B(x,4a) = X \cap B(x,4a)$
be given. Then 
$$
|\pi^\perp_x(y) - \pi^\perp_x(x)| = \dist(y,P_x)
\leq \varepsilon
\leqno (2.8)
$$
because $P_x$ goes through $x$, and then because $y \in X \i E$
and by (2.4) (also see the definition (1.2)).
Similarly, $|\pi^\perp_x(z) - \pi^\perp_x(x)| \leq \varepsilon$.
If $y = z$, the inequality in (2.7) is obvious.
Otherwise, $|y-z| \geq a$ by definition of $X$, and
$$
|\pi_x^\perp(y)-\pi_x^\perp(z)| \leq 2 \varepsilon
\leq 2 a^{-1}\varepsilon |y-z|;
\leqno (2.9)
$$
then $|\pi_x(y)-\pi_x(z)| \geq |y-z|
- |\pi_x^\perp(y)-\pi_x^\perp(z)| \geq |y-z|/2$
and (2.7) holds with $A'_0 = 4a^{-1} = 128$.
So there are Lipschitz functions 
$F_{0,x}$ such that (2.6) holds for $j=0$.

Next assume that $0 \leq j < N$ and that we constructed $S_j$
and the $F_{j,x}$ with the property (2.6). We take
$$
S_{j+1} = S_j \cup 
\bigcup_{w\in X_{j+1}} {\cal G}(F_{j,w}) \cap B(w,3a). 
\leqno (2.10)
$$
Let us now prove (2.7) for $j+1$. Let $x\in X$
and $y,z \in S_{j+1} \cap B(x,4a)$ be given; we need to
prove that
$$
|\pi_x^\perp(y)-\pi_x^\perp(z)| 
\leq A'_{j+1} \varepsilon |\pi_x(y)-\pi_x(z)|.
\leqno (2.11)
$$
If $y, z$ both lie on $S_j$, (2.11) simply follows 
from (2.7) for $j$. So we may assume that one of the
two points (say $y$ for definiteness) lies in 
${\cal G}(F_{j,w}) \cap B(w,3a)$ for some $w\in X_{j+1}$.

\ssi{\bf Case 1.}
We first assume that $z\in {\cal G}(F_{j,w})$.
Then
$$
|\pi_w^\perp(y)-\pi_w^\perp(z)| 
\leq A_j \varepsilon |\pi_w(y)-\pi_w(z)|,
\leqno (2.12)
$$
just because $F_{j,w}$ is $A_j \varepsilon$-Lipschitz.
This is not exactly (2.11), because we project on different
planes, but the difference will be small. 
Indeed, $|w-x| \leq |w-y|+|y-x| \leq 3a+4a < 1/2$, so (2.5)
says that
$$
d_{x,1/4}(P_x,P_w) \leq 10 \varepsilon.
\leqno (2.13)
$$ 
The reader will not be surprised to learn that
(2.11) follows from (2.12) and  (2.13), but let us check this
anyway. Denote by $\wt P_x$ and $\wt P_w$ the vector planes
parallel to $P_x$ and $P_w$, and by $\wt \pi_x$ 
(respectively $\wt \pi_w$) the orthogonal projection on $\wt P_x$ 
(respectively $\wt P_w$); we want to check that
$$
|| \wt\pi_x-\wt\pi_w|| \leq 100 \varepsilon.
\leqno (2.14)
$$ 
Let $v\in \R^n$ be given; we want to estimate
$| \wt\pi_x(v)-\wt\pi_w(v)|$, and we may as well 
assume that $v$ is a unit vector. We start with the
case when $v\in \wt P_x$, and we first check that
$$
\dist(v,\wt P_w) \leq 25 \varepsilon.
\leqno (2.15)
$$
Observe that both $x$ and $x+v/5$ lie in $P_x \cap B(x,1/4)$, 
so by (2.13) we can find $x'$ and $x''$ in $P_w$ such that
$|x'-x| \leq 10\varepsilon/4$ and $|x''-x-v/5| \leq 10\varepsilon/4$;
then $v'= 5(x''-x')$ lies in $\wt P_w$, and
$|v'-v| = 5|(x'' - x - v/5) + (x - x')| \leq 25 \varepsilon$, 
as needed for (2.15).

Because of (2.15), we get that
$|\wt\pi_x(v)-\wt\pi_w(v)| = |v-\wt\pi_w(v)| = \dist(v,\wt P_w)
\leq 25 \varepsilon$ for every unit vector $v\in \wt P_x$.

Now suppose $v$ is a unit vector in $P_x^\perp$. We set 
$\xi = \wt\pi_w(v)$, and check that 
$$
\dist(\xi,\wt P_x) \leq 25 \varepsilon. 
\leqno (2.16)
$$
Pick $x' \in P_w$ such that $|x'-x| \leq 10\varepsilon/4$ 
(use (2.13) as before), 
and, since $x'+ \xi/5 \in P_w \cap B(x,1/4)$, use (2.13)
to find $y\in P_x$ such that $|y-x'-\xi/5| \leq 10\varepsilon/4$.
Recall that $x\in P_x$; then $5(y-x) \in \wt P_x$, and 
$|5(y-x) - \xi| \leq 5[|y-x'-\xi/5|+ |x'-x|] \leq 25\varepsilon$,
which proves (2.16). Thus 
$$
|\xi - \wt \pi_x(\xi)| = \dist(\xi,\wt P_x)
\leq 25 \varepsilon. 
\leqno (2.17)
$$
Also, for each unit vector $e \in \wt P_x$, use (2.15) to find 
$e'\in \wt P_w$, with $|e'-e| \leq 25 \varepsilon$, 
and  observe that
$$\eqalign{
|\langle e, \xi \rangle| 
&\leq |\langle e', \xi \rangle| + |e'-e||\xi| 
\leq |\langle e', \xi \rangle| + 25 \varepsilon
= |\langle e', v \rangle| + 25 \varepsilon
\cr&
\leq |\langle e, v \rangle| + |e'-e| + 25 \varepsilon
\leq |\langle e, v \rangle| + 50 \varepsilon = 50 \varepsilon
}\leqno (2.18)
$$
because $\xi -v \perp e'$ by definition of $\xi=\wt \pi_w(v)$,
and $v\in P_x^\perp$. By (2.18), $|\wt \pi_x(\xi)| \leq 50 
\varepsilon$ and, adding up with (2.17),
$|\xi| \leq 75 \varepsilon$. Then $|\wt\pi_x(v)-\wt\pi_w(v)|
=|\wt \pi_w(v)| = |\xi| \leq 75 \varepsilon$
for every unit vector $v \in P_x^\perp$. Since we also
have $|\wt\pi_x(v)-\wt\pi_w(v)| \leq 25 \varepsilon$
for every unit vector $v \in \wt P_x$, we get that
$|\wt\pi_x(v)-\wt\pi_w(v)| \leq 100 \varepsilon |v|$
for every $v\in \R^n$. That is, (2.14) holds. 

Return to our $y$ and $z  \in S_{j+1} \cap B(x,4a)$, recall that 
$\pi_x^\perp$ denotes the orthogonal projection on 
the vector space $P_x^\perp$, define $\pi_w^\perp$ similarly,
and observe that
$$\leqalignno{
|(\pi_x^\perp(y)-\pi_x^\perp(z)) - \pi_w^\perp(y) -\pi_z^\perp(z))| 
&= |\pi_x^\perp(y-z) - \pi_w^\perp(y-z)|
\cr&
= |\wt \pi_x(y-z) - \wt\pi_w(y-z)|
\cr&
\leq || \wt\pi_x-\wt\pi_w|| \, |y-z|
\leq 100 \varepsilon |y-z|
& (2.19)
}
$$
because $\pi_x^\perp + \wt\pi_x = I = \pi_w^\perp + \wt\pi_w$.
Now (2.11) follows from (2.12) and (2.19),
as soon as we take $A'_{j+1} \geq A_j + 100$. This completes
our proof of (2.11) when $z\in {\cal G}(F_{j,w})$.

\ssi{\bf Case 2.}
We still need to check (2.11) when $z \notin {\cal G}(F_{j,w})$.
Notice that $z\in {\cal G}(F_{j,w'})\cap B(w',3a)$
for some other $w' \in X_{j+1}$ is impossible, because
we would get that 
$$
|w-w'| \leq |w-y|+|y-x|+|x-z|+|z-w'|
\leq 3a + 4a + 4a + 3a < 16a, 
\leqno (2.20)
$$
which is forbidden by (2.3).
By (2.10) and because $z\notin {\cal G}(F_{j,w})$,
$z \in S_j$. In addition, by (2.6) and 
because $z\notin {\cal G}(F_{j,w})$, $z$ lies
out of $B(w,4a)$, and so $|z-y| \geq |z-w|-|w-y| \geq
4a-|w-y| \geq a$ (recall from the line below (2.11) that 
$y\in {\cal G}(F_{j,w}) \cap B(w,3a)$). If we prove that
$$
|\pi_x^\perp(y)-\pi_x^\perp(z)| \leq C \varepsilon
\leqno (2.21)
$$
for some $C$ that depends on $n$ and $A'_j$,
(2.11) will follow as soon because (2.21) also implies that
$|\pi_x(y)-\pi_x(z)| \geq a/2$.
Observe that
$$
|\pi_x^\perp(x)-\pi_x^\perp(z)|
\leq 4 a A'_j \varepsilon
\leqno (2.22)
$$
by (2.7), and because $z\in B(x,4a)$ by definition and
we just checked that $z\in S_j$. We are left with
$$\eqalign{
|\pi_x^\perp(y)-\pi_x^\perp(x)|
&\leq |\pi_x^\perp(y)-\pi_x^\perp(w)|+|\pi_x^\perp(w)-\pi_x^\perp(x)|
\cr&
= |\pi_x^\perp(y)-\pi_x^\perp(w)|+\dist(w,P_x)
\cr&
\leq |\pi_x^\perp(y)-\pi_x^\perp(w)|+\varepsilon
}\leqno (2.23)
$$
by (2.4) for $x$, and because $|w-x| \leq |w-y|+|y-x| \leq 7a < 1/2$. 
Let us again use (2.14) (its proof is valid also in case 2; 
in fact neither the statement nor the proof involves $z$); 
this yields
$$\eqalign{
|\pi_x^\perp(y)-\pi_x^\perp(w)| 
&= |\pi_x^\perp(y-w)| 
\leq |\pi_w^\perp(y-w)| + |y-w|\,||\pi_x^\perp - \pi_w^\perp||
\cr&
= |\pi_w^\perp(y-w)| + |y-w|\,||\pi_x - \pi_w||
\cr&
\leq |\pi_w^\perp(y-w)| + 100 \varepsilon |y-w|
}\leqno (2.24)
$$
because $\pi_x^\perp$ is a linear projection and
$\pi_x^\perp + \wt\pi_x = I = \pi_w^\perp + \wt\pi_w$.
Now recall that $y\in {\cal G}(F_{j,w}) \cap B(w,3a)$
(see below (2.11)), and so does $w$ 
(because $w \in X = S_0 \i S_j$, and by (2.6)); then
$$\eqalign{
|\pi_w^\perp(y-w)| &=  |\pi_w^\perp(y)-\pi_w^\perp(w)|
= |F_{j,w}(\pi_w(y))-F_{j,w}(\pi_w(w))|
\cr&
\leq A_j \varepsilon |\pi_w(y)-\pi_w(w)|
\leq A_j \varepsilon |y-w| \leq A_j \varepsilon
}\leqno (2.25)
$$
because $F_{j,w}$ is $A_j \varepsilon$-Lipschitz.
Altogether, $|\pi_x^\perp(y)-\pi_x^\perp(z)| \leq
4 a A'_j \varepsilon + \varepsilon + 300a \varepsilon + A_j 
\varepsilon$, by (2.22)-(2.25), which proves (2.21) and
then (2.11) in our second case.

This completes our proof of (2.7) for $j+1$.
Thus we have the local Lipschitz description of $S_{j+1}$,
with (2.6), and this completes our construction of 
sets $S_j$ by induction.

\ms
Let us now check that for $1 \leq j \leq N$ and $x\in X_j$,
$$
S_N \cap B(x,2a) = {\cal G}(F_{j-1,x}) \cap B(x,2a).
\leqno (2.26)
$$
Set ${\cal G} = {\cal G}(F_{j-1,x})$.
By (2.10), ${\cal G} \cap B(x,3a) \i S_j \i S_N$
so we get a first inclusion. Notice also that since the
Lipschitz constant is small and $\cal G$ goes through $x$,
$$
\pi_x({\cal G} \cap B(x,3a))
\hbox{ contains } P_x \cap B(x,2a). 
\leqno (2.27)
$$
Now let $z\in S_N \cap B(x,2a)$ be given. 
Then $\pi_x(z) \in P_x \cap B(x,2a)$; by (2.27) we can find
$w\in {\cal G} \cap B(x,3a)$ such that $\pi_x(w) = \pi_x(z)$. 
Then $w\in S_N \cap B(x,3a)$.
But (2.6) says that $S_N \cap B(x,4a)$
is contained in a Lipschitz graph, so $\pi_x$ is injective on
$S_N \cap B(x,4a)$, hence $w=z$ and $z\in {\cal G}$.
So $S_N \cap B(x,2a) \i {\cal G}$, as needed for (2.26).

\ms
Now we check the properties (1.4) and (1.5) for
$S_N$. First let $y\in E$ be given. By (2.1), we
can find $x\in X$ such that $|x-y| \leq a$. Let
$j$ be such that $x\in X_j$, and set
$w = \pi_x(y)$ and $z = w + F_{j-1,x}(w) \in {\cal G}(F_{j-1,x})$.
We know from (2.26) that $z\in S_N$, and so
$$\eqalign{
\dist(y,S_N) &\leq |y-z| = |\pi_x^\perp(y)-\pi_x^\perp(z)|
= |\pi_x^\perp(y)-F_{j-1,x}(w)|
\cr&
\leq |\pi_x^\perp(y)-\pi_x^\perp(x)| + |\pi_x^\perp(x)-F_{j-1,x}(w)|
\cr&
= \dist(y,P_x) + |\pi_x^\perp(x)-F_{j-1,x}(w)|
}\leqno (2.28)
$$
because $\pi_x(y) = \pi_x(z) = w$ and then $x\in P_x$.
But $x \in X \i S_N$, so $x\in {\cal G}(F_{j-1,x})$
by (2.26) and hence $F_{j-1,x}(x) = \pi_x^\perp(x)$. 
Thus
$$
|\pi_x^\perp(x)-F_{j-1,x}(w)| = |F_{j-1,x}(x)-F_{j-1,x}(w)|
\leq A_{j-1} \varepsilon |x-w| \leq C \varepsilon
\leqno (2.29)
$$
because $F_{j-1,x}$ is $A_{j-1} \varepsilon$-Lipschitz.
Since $\dist(y,P_x) \leq \varepsilon$ by (2.4),
we deduce from (2.28) that $\dist(y,S_N) \leq C\varepsilon$.
So (1.4) holds.

For (1.5), let $y\in S_N$ be given, and let 
$j(y)$ be the first index such that $y\in S_j$; 
if $j(y)=0$, $y$ lies in $X \i E$, and $\dist(y,E) = 0$.
Otherwise, (2.10) says that 
$y \in {\cal G}(F_{j-1,x}) \cap B(x,3a)$ for some $x\in X_j$. 
Set $w = \pi_x(y)$; then $w \in P_x \cap B(x,3a)$ and hence
$\dist(w,E) \leq \varepsilon$ by (2.4). Now
$$\eqalign{
\dist(y,E) 
&\leq |y-w| + \dist(w,E) 
\leq |y-w| + \varepsilon
= |\pi_x^\perp(y-w)|
= |\pi_x^\perp(y-x)|
\cr&
= |F_{j-1,x}(w)-F_{j-1,x}(x)| + \varepsilon
\leq A_{j-1} \varepsilon |x-w| +  \varepsilon \leq C \varepsilon,
}\leqno (2.30)
$$
where we used the facts that $w = \pi_x(y)$, that
$\pi_x^\perp(w-x)=0$ because $x$ and $w$ lie in $P_x$, that 
$y = w + F_{j-1,x}(w)$ and $x = x+F_{j-1,x}(w)$
(both points lie in ${\cal G}(F_{j-1,x})$; also see the 
definition (1.7)), and that $F_{j-1,x}$ is 
$A_{j-1} \varepsilon$-Lipschitz. This proves (1.5) for $S_N$.

The local Lipschitz description (1.7)-(1.8) for $S_N$ 
(but with no extra smoothness yet) will easily follow
from (2.26). Indeed, if $y\in S_N$, (1.5) and (2.1) imply
that $\dist(y,X) \leq a+\varepsilon$; we choose $j\in [1,N]$
and $x\in X_j$ so that $|x-y| \leq a + \varepsilon$ 
and try to take $P_y = P_x$ and $F_y = F_{j-1,x}$; then (1.8)
follows from (2.26) because $B(y,10^{-2}) \i B(x,2a)$.
(Recall that $a = {1 \over 32}$.)
We also get (1.9) for $k=1$
(i.e., $F_y = F_{j-1,x}$ is $C \varepsilon$-Lipschitz).
Now $P_x$ does not necessarily pass through $y$, but 
we can just translate it, modify $F_y$ accordingly, and
get an acceptable pair $(P_y,F_y)$ that works.

\ms
Our proof is not complete yet, because we want our smooth 
surface $\Sigma_0$ to satisfy (1.9) for all $k$.
If we were only interested in a finite number of
derivatives, we could modify the argument above, use
Whitney extensions with functions of class $C^k$, and
conclude as before.
But since we decided to require an infinite
number of derivatives, the simplest seems to
start from $S_N$ and smooth it out. Notice that
if we do not move it by more than $C\varepsilon$,
(1.4) and (1.5) stay true (with larger constants),
so the reader should not expect trouble here.

Let us rapidly describe the argument. The details would
be standard and boring, so we skip them.
We shall use the same sets $X_j$, $1 \leq j \leq N$,
as before, and define sets 
$T_0 = S_N, T_1, \ldots$, and finally $\Sigma_0 = T_N$, by induction. 
For $j = 1, \ldots N$, we obtain $T_j$ from $T_{j-1}$
by smoothing $T_j$ in the balls $B(x,{3a \over 2})$, $x\in X_j$.
Since these balls are far from each other, we can do this
independently in all these balls. We start from a description of 
$T_{j-1}$ as a Lipschitz graph near $B(x,{3a \over 2})$
(which we get from an induction assumption), use a convolution with
a smooth function with a support of small diameter to smooth out
$T_{j-1}$ in $B(x,{3a \over 2})$, interpolate smoothly in a small ring
near $\partial B(x,{3a \over 2})$, and get the next $T_j$.
Each time we get a Lipschitz graph with a slightly worse
constant, but this does not matter. The smoothness that
has been produced in the previous stages is preserved too, 
and at the end, $\Sigma_0 = T_N$ is the desired smooth set.
This completes our proof of Theorem 1.10.
\qed

\msi
{\bf Remark 2.31.}
In [DS], one also obtains parameterizations 
of sets that contain $E$, assuming that for $x\in E$
and $0 < r \leq 1$, we can find a $d$-dimensional affine
plane $P(x,r)$ through $x$ such that
$$
\dist(y,P(x,r)) \leq \varepsilon r
\ \hbox{ for } y \in E\cap B(x,r),
\leqno (2.32)
$$
and also that the $P(x,r)$ vary sufficiently slowly. 
Here we can do something similar, suppose that for $x\in E$, 
we can find a $d$-plane $P_x$ through $x$, such that
$$
\dist(y,P_x) \leq \varepsilon 
\ \hbox{ for } y \in E\cap B(x,1),
\leqno (2.33)
$$
and such that (2.5) holds. Then our construction gives a set
$S_N$, with a Lipschitz description near $E$, as in (2.26)
included. We also get (1.4) as before, but of course not (1.5).
But slightly more unpleasantly, $S_N$ has a boundary; that is,
we only get (2.26) for balls centered on $X \i E$, but not
necessarily for balls centered on points of $S_N$ far from $E$.
For instance, the parts near $E$ may be far from each other
and then $S_N$ stops somewhere in between. It is possible 
that we can connect all the boundary pieces of $S_N$, and get a 
larger set $\Sigma_0$ which is smooth and without boundary, but
this could be unpleasant to do. The author did not find any 
obvious topological obstruction either.

\msi
{\bf Remark 2.34.}
When $E$ is unbounded, we can probably prove a variant of 
Theorem 1.10, where the radius $r_0$ depends slowly on the
location of the balls. We shall not give any detail.

\medskip\noindent
{\bf 3. Approximation from the inside of Reifenberg-flat domains.}
\ms

In this section we apply Theorem 1.10 to prove Theorem 1.18. 
Let $E$ be as in the statement. For each $r \in (0,r_0]$,
we can apply Theorem 1.10 (with $r_0$ replaced with $r$),
and we get a smooth surface $\Sigma_{r}$ that satisfies
the (1.4)-(1.9) with the radius $r$. The general plan is
easy to guess: we shall move $\Sigma_{r}$ in the direction of
the unit normal to get two new surfaces 
$\Sigma_{r,1}$ and $\Sigma_{r,2}$ (one on each side
of $\Sigma_{r}$), which will bound the desired domains 
$W_{r,1}$ and $W_{r,2}$. The slightly unpleasant part of the
argument will be to take care of the topology. In particular,
we shall give a more direct proof of (1.20), that relies only
on the simple argument of [S] (rather that the full Reifenberg 
machine), and which is a little reminiscent of arguments in 
Chapter II.4 of [DS]. 

Let us first check that for $0 < r \leq r_0$,
$$
\Sigma_{r} \hbox{ is connected.}
\leqno (3.1)
$$
When $E$ is connected, we can use the proof of (1.25)
to get this directly, but let us also prove (3.1) under
the weaker assumption that $E$ is ${r_0 \over 20}$-connected 
(as in Remark 1.26). The proof of (1.25) shows
that $\Sigma_{r_0}$ is connected, and then (1.4) and (1.5)
imply that $E$ is $10 C_0 \varepsilon r_0$-connected.
Then set $r_1 = r_0/10$, use the fact that 
(if $\varepsilon < \varepsilon_2$ is small enough) $E$ is
${r_1 \over 20}$-connected to prove that $\Sigma_{r_1}$ is 
connected, and continue. We get that $E$ is $\rho$-connected
for every $\rho > 0$, and it is easy to deduce from this
that it is connected. So the weaker ${r_0 \over 20}$-connectedness
implies that $E$ is connected, and we get a full proof of (3.1).

By [S], $\Sigma_{r}$ is orientable, which means that there is a 
continuous (in fact, smooth) choice of unit normal to $\Sigma_{r}$;
denote by $n(y)$ this unit normal, computed at the point $y\in \Sigma_r$.
The two surfaces that we want to consider are 
$$
\Sigma_{r,1} = \big\{ z \in \R^n \,; \, z = y + 5C_0\varepsilon r n(y)
\hbox{ for some } y\in \Sigma_{r} \big\}
\leqno (3.2)
$$
and
$$
\Sigma_{r,2} = \big\{ z \in \R^n \,; \, z = y - 5C_0\varepsilon r n(y)
\hbox{ for some } y\in \Sigma_{r} \big\},
\leqno (3.3)
$$
where $C_0$ is as in (1.4) and (1.5).
These will eventually be our $\d W_{r,j}$, but we shall take care
of the topology  later.
Our first task is to show that the $\Sigma_{r,j}$ satisfy the constraints 
(1.22)-(1.24) that we require from the boundaries $\d W_{r,j}$.
For this we will use the good description of $\Sigma_r$ that 
comes from (1.6)-(1.9).

Let $y \in \Sigma_r$ be given; we want to study the $\Sigma_{r,j}$
near $y$. Recall from (1.6)-(1.9) that there is a hyperplane $P_y$ through $y$
and a smooth function $F_y : P_y \to P_y^\perp$, such that
$\Sigma_{r}$ coincides with the graph of $F_y$ in $B(y,10^{-2} r)$.
Let $\cal G$ denote the graph of $F_y$ over $P_y$, defined as
in (1.7). We want to compute a unit normal to $\cal G$, and this
will be simpler after a change of basis. Without loss of generality,
we can assume that $y=0$, and that we chose an orthonormal basis 
$(e_1, \ldots, e_n)$ of $\R^n$ such that $P_y$ 
is the plane of equation $x_n = 0$. Then, for $w\in P_y$,
a unit normal to ${\cal G}$ at $w+F_y(w)$ is 
$$
\nu(w) = a(w) \big(e_n - \sum_{j=1}^{n-1} 
e_j{\partial F_y \over \partial w_j}(w) \big), \hbox{ with }
a(w) = \big|e_n - \sum_{j=1}^{n-1} 
e_j {\partial F_y \over \partial w_j}(w) \big|^{-1}.
\leqno (3.4)
$$
Let us first study $\cal G$ and the analogues of the $\Sigma_{r,j}$
for $\cal G$, and then we will use (1.7) to return to the $\Sigma_{r,j}$.
For $w \in P_y$, set
$$
z(w) = w+F_y(w) \in {\cal G},
\leqno (3.5)
$$
$$
\xi_\pm(w) = z(w) \pm 5 C_0 \varepsilon r \nu(w)
= w + F_y(w) \pm 5 C_0 \varepsilon r \nu(w),
\leqno (3.6)
$$
and then 
$$
{\cal G}_\pm = \big\{ \xi_{\pm}(w)  \,; \, w \in P_y \big\}.
\leqno (3.7)
$$
We want to show that ${\cal G}_\pm$ is 
a smooth Lipschitz graph, and for this we will just need to estimate 
some derivatives and apply the implicit function theorem.

Recall from (1.9) that $F_y$ is $\Lambda_1\varepsilon$-Lipschitz,
so $\big|\sum_{j=1}^{n-1} e_j {\partial F_y \over \partial w_j}(w) \big| 
\leq \Lambda_1\varepsilon$; then the normalizing factor $a(w)$ in 
(3.4) stays close to $1$ (if $\varepsilon_2$ is chosen small enough),
$$
|\nu(w) - e_n| \leq C \varepsilon,
\leqno (3.8)
$$
and (1.9) and (3.4) yield
$$
|D^k \nu(w)| \leq C_k \varepsilon r^{-k}
\leqno (3.9)
$$
for $w \in P_y$ and $k \geq 1$. Here and below, $C_k$
denotes a constant that depends only on $k$ and $n$.
From (3.6), (3.9), and (1.9) we also deduce that
$$
|D^k [\xi_\pm(w)-w]| \leq C_k \varepsilon r^{1-k}
\leqno (3.10)
$$
(where for $k=0$, we also use the fact that $F_y(y)=0$).
Let $\pi$ denote the orthogonal projection on $P_y$,
and set $\psi_\pm = \pi \circ \xi_\pm$; then
$$
|D\psi_\pm(w) - I| \leq |D\pi \circ DF_y(w)|
+ 5 C_0 \varepsilon r |D\pi \circ D\nu(w)|
\leq C \varepsilon
\leqno (3.11)
$$
by (3.6), (1.9), and (3.9). 
By the implicit function theorem (or rather the contracting fixpoint 
theorem), $\psi_\pm$ is a bijection from $P_y$ to $P_y$; denote by 
$\varphi_\pm : P_y \to P_y$ its inverse. Notice that 
$\varphi_\pm$ is differentiable and
$$
D\varphi_\pm(u) = D\psi_\pm^{-1}(\varphi_\pm(u));
\leqno (3.12)
$$
then $\varphi_\pm$ is $(1+C\varepsilon)$-Lipschitz, by (3.11).
In addition, (3.10) also holds for $\psi_\pm$ (trivially),
and then for its inverse $\varphi_\pm$ (because of (3.12)),
and also the composition $f_{\pm} = \xi_\pm  \circ \varphi_\pm$.
That is, 
$$
||D^k [f_\pm-I]||_\infty \leq C_k \varepsilon r^{1-k}.
\leqno (3.13)
$$
It is easy to see that 
$$
{\cal G}_\pm \hbox{ is the graph of }
f_{\pm} - I : P_y \to P_y^\perp, 
\leqno (3.14)
$$
so we have a good Lipschitz graph description of ${\cal G}_\pm$, 
of the same type as in (1.6)-(1.9) and (1.23)-(1.24).
Let us also check that
$$
4 C_0 \varepsilon r  \leq \dist(\xi,{\cal G})
\leq 5 C_0 \varepsilon r 
\ \hbox{ for } \xi \in {\cal G}_\pm.
\leqno (3.15)
$$
Let $\xi \in {\cal G}_\pm$ be given, and write 
$\xi =\xi_\pm(w) = z(w) \pm 5 C_0 \varepsilon r \nu(w)$
for some $w \in P_y$. The second inequality is trivial
because $z(w) \in {\cal G}$; for the first one
we shall just use the last coordinates. We proceed
by contradiction and suppose that $B(\xi, 4C_0 \varepsilon r)$ 
contains a point $z_1 \in {\cal G}$.
Write $z_1 = z(w_1)$ and observe that
$$\eqalign{
|\xi - z_1| &\geq |\langle \xi - z_1, e_n \rangle|
\geq |\langle \xi - z(w), e_n \rangle| - |\langle z(w)-z_1, e_n \rangle|
\cr&
= |\langle 5 C_0 \varepsilon r \nu(w), e_n \rangle| 
- |\langle z(w)-z(w_1), e_n \rangle|
\cr&
\geq 5 C_0 \varepsilon r
- 5 C_0 \varepsilon r |\nu(w) - e_n| - |F_y(w)-F_y(w_1)|
\cr&
\geq 5 C_0 \varepsilon r - C \varepsilon^2 r - \Lambda_1 \varepsilon |w-w_1|
\geq  5 C_0 \varepsilon r - C \varepsilon^2 r
\geq 4 C_0 \varepsilon r
}\leqno (3.16)
$$
by (3.6), (3.8), because $F_y$ is $\Lambda_1 \varepsilon$-Lipschitz and 
$|w-w_1| \leq  |z(w)-z(w_1)| \leq |z(w)-\xi| + |\xi-z_1|
\leq 10C_0 \varepsilon r$, and if $\varepsilon_2$ is small enough. 
This contradiction completes our proof of (3.15).

\ms
Now we compare the ${\cal G}_\pm$ and the $\Sigma_{r,j}$.
Since by (1.7), $\Sigma_{r} = {\cal G}$ inside of $B(y,10^{-2} r)$,
$z(w) = w+F_y(w)$ is also a parameterization of $\Sigma_r$ near $y$, and
$\nu(w) = \pm n(z(w))$ in a neighborhood of $y$. Possibly at the price
of replacing $e_n$  with $-e_n$ and $F_y$ by $-F_y$, we may assume 
that the sign is $+$. Then (1.7) and a small connectedness argument 
show that
$$
\nu(w) = n(z(w))
\ \hbox{ for } w \in P_y \cap B(y,10^{-2} r/2).
\leqno (3.17)
$$
Let us check that
$$
\Sigma_{r,1} \cap B(y,10^{-2}r/3) = {\cal G}_+ \cap B(y,10^{-2}r/3).
\leqno (3.18)
$$
If $\xi \in {\cal G}_+ \cap B(y,10^{-2}r/3)$, then by (3.7)
$\xi = \xi_+(w) = z(w) + 5 C_0 \varepsilon r \nu(w)$ 
for some $w \in P_y$. Notice that $w = \pi(z(w)) \in B(y,10^{-2}r/2)$
because $|\xi-z(w)| = 5 C_0 \varepsilon r$, so $\nu(w) = n(z(w))$
by (3.17). Also, $z(w) \in {\cal G}_+ \cap B(y,10^{-2}r/2) \i 
\Sigma_r$, by (1.8), and now 
$\xi = z(w) + 5 C_0 \varepsilon r n(z(w))$ 
lies in $\Sigma_{r,1}$, by (3.2).

Conversely, if $\xi \in \Sigma_{r,1} \cap B(y,10^{-2}r/3)$,
then we can write $\xi = z + 5 C_0 \varepsilon r n(z)$
for some $z \in \Sigma_{r}$. Obviously $z\in B(y,10^{-2}r/2)$,
so $\Sigma_{r} = {\cal G}$ near $z$ (by (1.8)); then 
$z = z(w)$ for some $w\in P_y$. In addition, 
$w = \pi(z) \in B(y,10^{-2}r/2)$, so $n(z) = \nu(w)$ by (3.7).
Thus $\xi = z(w) + 5 C_0 \varepsilon r \nu(w) \in {\cal G}_+$, 
as needed for (3.18). The proof of (3.18) also yields
$$
\Sigma_{r,2} \cap B(y,10^{-2}r/3) = {\cal G}_- \cap B(y,10^{-2}r/3).
\leqno (3.19)
$$
Let us now deduce from this that
$$
4 C_0 \varepsilon r \leq \dist(\xi,\Sigma_{r}) \leq 5 C_0 \varepsilon r
\ \hbox{ for } \xi \in \Sigma_{r,1} \cup \Sigma_{r,2}.
\leqno (3.20)
$$
Let $\xi$ lie in some $\Sigma_{r,j}$, and use the definition (3.2) or 
(3.3) to write $\xi = y \pm 5 C_0 \varepsilon r n(y)$ for some 
$y\in \Sigma_{r}$; then $\xi \in B(y,10^{-2}r/3)$ and, with the notation above, 
$\xi \in {\cal G}_+ \cup {\cal G}_-$ (by (3.18) or (3.19)); 
now (3.20) follows from (3.15) and the fact that $\Sigma_r = {\cal G}$
in $B(y,10^{-2}r)$.

\ms 
Let us start to deal with the topology. 
Recall from the discussion below (3.1) that $\Sigma_r$ is connected, 
hence orientable. In addition,
$$
\R^n \sm \Sigma_{r} \hbox{ has exactly two connected components;}
\leqno (3.21)
$$
see for instance [S]. For each $y\in \Sigma_r$, we know from
(1.8) that in $B(y,10^{-2}r)$, $\Sigma_r$ coincides with a
Lipschitz graph ${\cal G} = {\cal G}(y)$. Denote by $U_1(y)$
the connected component of $\R^n \sm \Sigma_{r}$ that contains
the points of $B(y,10^{-3}r) \sm \Sigma_r$ that lie in the
direction of $n(y)$; in the coordinates that we used to
describe $\Sigma_r = {\cal G}$ near $y$, this is the 
component that contains the part of $B(y,10^{-3}r) \sm {\cal G}$
above ${\cal G}$. It is clear that $U_1(y)$ is locally constant
on $\Sigma_r$, i.e., that $U_1(y)=U_1(y')$ when $|y'-y| \leq 10^{-3}r$.
Since $\Sigma_r$ is connected, $U_1(y)$ is constant; we shall denote 
by $U_{r,1}$ its constant value on $\Sigma_r$.
Similarly, we can define $U_{r,2}$ to be the connected component
of $\R^n \sm \Sigma_{r}$ that contains the points of 
$B(y,10^{-3}r) \sm \Sigma_r$ that lie in the direction of $-n(y)$,
because this component does not depend on $y$.

Each of the two components of $\R^n \sm \Sigma_{r}$ gets arbitrarily close to
$\Sigma_r$, hence it meets some $B(y,10^{-3}r)$ and it is equal to 
$U_{r,1}$ or $U_{r,2}$. This proves that $U_{r,1} \neq U_{r,2}$; and
these are the two components of $\R^n \sm \Sigma_{r}$.

We have a local Lipschitz description of the two $\Sigma_{r,j}$,
that comes from (3.18) and (3.19) (also see (3.13)-(3.14)), 
and from the connectedness of $\Sigma_r$ we easily deduce 
that $\Sigma_{r,j}$ is connected: given $z_1, z_2 \in \Sigma_{r,j}$, 
pick $y_i \in \Sigma_r \cap B(z_i,10^{-3}r)$, connect $y_1$ to $y_2$
by a $(10^{-3}r)$-chain of points $y\in \Sigma_r$, use (3.18) or (3.19)
to connect $z_1$ to successive points $z\in \Sigma_{r,j} \cap 
B(z,10^{-3}r)$ by paths in $\Sigma_{r,j}$, and finally
connect the last $z$ to $z_2$. 

By [S], $\R^n \sm \Sigma_{r,j}$ has exactly two 
connected components, which we denote by $W_{r,j}$ and $W'_{r,j}$. 
Since $\Sigma_r$ is connected and does not meet $\Sigma_{r,j}$
(by 3.20)), it is contained in one of these components, and
we choose their names so that 
$$
\Sigma_r \i W'_{r,j}
\ \hbox{ and } \Sigma_r \cap W_{r,j} = \emptyset.
\leqno (3.22)
$$
For $y\in \Sigma_r$, we can use (3.18) and (3.19) to localize these 
components near $y$. We have three Lipschitz graphs,
${\cal G}_-$, ${\cal G}$, and ${\cal G}_+$, which do not meet
(by (3.15)), and ${\cal G}_-$ and ${\cal G}_+$ cut $\R^n$ into
three open region: the region $R_-$ below ${\cal G}_-$,
the region $R_0$ between ${\cal G}_-$ and ${\cal G}_+$
(and which contains ${\cal G}$), and the region $R_+$ above 
${\cal G}_+$. Let us check that
$$
\Sigma_{r,j} \cap B(y,10^{-2}r/4) \i U_{r,j}
\hbox{ for } j = 1,2,
\leqno (3.23)
$$
$$
R_0 \cap B(y,10^{-2}r/4) \i W'_{r,1} \cup W'_{r,2},
\leqno (3.24)
$$
$$
R_+ \cap B(y,10^{-2}r/4) \i W_{r,1}, 
\leqno (3.25)
$$
and 
$$
R_- \cap B(y,10^{-2}r/4) \i W_{r,2}.
\leqno (3.26)
$$
Let us first prove (3.23) for $j=1$. Recall from
the definitions below (3.21) that $U_{r,1}$ is the common value, 
for $y\in \Sigma_r$, of $U_{1}(y)$, which is the component of 
$\R^n \sm \Sigma_r$ that contains the points of 
$B(y,10^{-3}r) \sm \Sigma_r$ that lie in the direction of $n(y)$.
Clearly this component also contains ${\cal G}_+ \cap 
B(y,10^{-2}r/4)$, which by (3.18) is the same as 
$\Sigma_{r,1} \cap B(y,10^{-2}r/4)$. This proves (3.23)
for $j=1$; the proof for $j=2$ is the same.

Next, if $x\in R_0 \cap B(y,10^{-2}r/4)$, the vertical line segment from
$x$ to the point of $\cal G$ with the same projection is 
contained in $B(y,10^{-2}r/3) \sm [{\cal G}_+ \cup {\cal G}_-]$,
hence (by (3.18) and (3.19)) in 
$B(y,10^{-2}r/3) \sm [\Sigma_{r,1} \cup \Sigma_{r,2}]$.
Since $y\in \Sigma_r \i W'_{r,j}$ for $j=1, 2$, we get that
$x\in W'_{r,j}$, which proves (3.24).

Let now check (3.25), starting with a specific $y$. 
Since $W_{r,1}$ is a component of $\R^n \sm \Sigma_{r,1}$,
we can find $\xi \in \Sigma_{r,1}$ such that 
$\dist(\xi, W_{r,1}) \leq C_0 \varepsilon r$.
Then (by (3.20)) we can pick $y \in \Sigma_r$
such that $|\xi-y| \leq 5 C_0 \varepsilon r$.
For this specific $y$, $\xi \in {\cal G}_+ \cap B(y,10^{-2}r/4)$
(by (3.18)), so 
$$
\dist({\cal G}_+ \cap B(y,10^{-2}r/4), W_{r,1}) \leq C_0 \varepsilon r. 
\leqno (3.27)
$$
Let $a \in W_{r,1}$ be such that 
$\dist(a,{\cal G}_+ \cap B(y,10^{-2}r/4) \leq 2 C_0 \varepsilon r$.
Observe that $a \notin \overline R_-$, because 
$\dist(\overline R_-,{\cal G}_+) \geq \dist({\cal G}_-,{\cal G}_+)
\geq \dist({\cal G}_-,{\cal G}) \geq 4C_0 \varepsilon r$
by (3.15).
But $a$ does not lie in $R_0$ (by (3.24)), nor in $\Sigma_{r,1}$
(by definition of $W_{r,1}$), so we get that $a\in R_+$.
It is then easy to connect $a$ to any point of 
$R_+ \cap B(y,10^{-2}r/4)$, by a path contained
in $B(y,10^{-2}r/3) \sm {\cal G}_+ = B(y,10^{-2}r/3) \sm \Sigma_{r,1}$
(by (3.18)). So (3.25) holds for this specific $y$.

But then every point of $\Sigma_{r,1} \cap B(y,10^{-2}r/4)$
lies very close to $W_{r,1}$, and (3.27) also holds for every 
$y'\in \Sigma_r$ such that $|y'-y| \leq 10^{-3}r$.
By the argument above, we also get (3.25) for every such $y'$
and eventually, because $\Sigma_r$ is connected, 
for every $y'\in \Sigma_r$. This proves (3.25) for all $y\in \Sigma_r$. 
The proof of (3.26) is the same.

We now have a definition of $W_{r,1}$ and $W_{r,2}$, and we can
start checking the conclusions of Theorem 1.18. First observe that
$$
\d W_{r,j} = \Sigma_{r,j} \, ;
\leqno (3.28)
$$
the first inclusion is obvious, and conversely every point of $\Sigma_{r,j}$
lies in $\d W_{r,j}$: we choose $y \in \Sigma_r$ not too far and use 
(3.25) or (3.26). From (3.20),
(1.4), and (1.5), we easily deduce that
$$
3 C_0 \varepsilon r \leq \dist(\xi,E) \leq 6 C_0 \varepsilon r
\ \hbox{ for } \xi \in \Sigma_{r,1} \cup \Sigma_{r,2},
\leqno (3.29)
$$
which proves (1.22) with $C_2 = 3 C_0$. For (1.21), we 
just need to check that 
$$
\big\{ x\in \R^n \, ; \, \dist(x,E) > 6 C_0 \varepsilon r \big\}
\i W_{r,1} \cup W_{r,2}.
\leqno (3.30)
$$
Suppose $\dist(x,E) > 6 C_0 \varepsilon r$, 
pick any point $a \in E$, and run along the segment $[x,a]$
until the first point $z$ such that $\dist(z,E) = 6C_0 \varepsilon r$. 
Thus 
$$
\dist(\xi,E) > 6C_0 \varepsilon r
\ \hbox{ for } \xi \in [x,z)
\leqno (3.31)
$$
and, by (3.29), $[x,z)$ does not meet $\Sigma_{r,1} \cup \Sigma_{r,2}$.
Pick $x'\in [x,z)$ close to $z$, so that
$\dist(x',E) \leq 7C_0 \varepsilon r$ and hence
$\dist(x',\Sigma_r) \leq 8C_0 \varepsilon r$ (by (1.4)).
Pick $y\in \Sigma_r$ such that $|y-x'| \leq 8C_0 \varepsilon r$,
and use $y$ to localize. If $x$ does not lie in $W_{r,1} \cup W_{r,2}$, 
then neither does $x'$, and since $x'\in B(y,10^{-3}r)$, 
(3.25) and (3.26) imply that $x'\in \overline R_0$. But (3.15) 
then implies that $\dist(x',{\cal G}) \leq 5C_0 \varepsilon r$.
In addition, ${\cal G} = \Sigma_r$ in $B(y,10^{-2} r)$ (by (1.8)), so
$\dist(x',\Sigma_r) \leq 5C_0 \varepsilon r$, and (by (1.5))
$\dist(x',E) \leq 6C_0 \varepsilon r$. This contradicts (3.31);
(3.30) and (1.21) follow.

Let us also check that
$$
\big\{ x\in \R^n \, ; \, \dist(x,E) < 3C_0 \varepsilon r \big\}
\i W'_{r,1} \cap W'_{r,2}
\i \R^n \sm [W_{r,1} \cup W_{r,2}].
\leqno (3.32)
$$
Suppose $\dist(x,E) < 3C_0 \varepsilon r$, and let $a\in E$
be such that $|x-a| < 3C_0 \varepsilon r$. 
Observe that $[x,a]$ does not meet $\Sigma_1 \cup \Sigma_2$, by (3.29).
Then use (1.4) to find $b\in \Sigma_r$ such that $|b-a| \leq C_0 
\varepsilon r$; then $[a,b]$ does not meet $\Sigma_1 \cup \Sigma_2$
either, again by (3.29). But $b\in W'_{r,1} \cap W'_{r,2}$, by (3.22),
so $x\in W'_{r,1} \cap W'_{r,2}$ too. The second part of (3.32) is 
immediate.

Now we establish the Lipschitz description of the $\d W_{r,j}$
in (1.23) and (1.24).
By (3.28), $\d W_{r,j} = \Sigma_{r,j}$.
Since $\Sigma_{r,j}$ stays close to $\Sigma_r$, by (3.20),
a good Lipschitz description of $\Sigma_{r,j}$ in each 
$B(y, 10^{-2} r/3)$ will give a good description in each 
$B(z, 10^{-3} r)$, $z\in \Sigma_{r,j}$. But in 
$B(y, 10^{-2} r/3)$, (3.18) or (3.19) says that
$\Sigma_{r,j}$ coincides with a ${\cal G}_\pm$, and now the
desired Lipschitz description follows from (3.13) and (3.14).

The fact that for $z\in \Sigma_{r,j}$,
$W_{r,j} \cap B(z,10^{-3}r)$ lies on only one side of $\Sigma_{r,j}$
can also be seen locally, in some $B(y, 10^{-2} r/3)$, $y\in \Sigma_r$, 
where it follows from (3.24)-(3.26); the proof is the same as for 
(3.28).

Our last tasks will be to show that $\R^n \sm E$ has 
exactly two connected components (as in (1.20)), and that 
for each $r$,
$$
\hbox{the two $W_{r,j}$ are contained in different components of 
$\R^n \sm E$.}
\leqno (3.33)
$$
Once we do this, Theorem 1.18 will follow, just after 
renaming the $W_{r,j}$ so that $W_{r,1}$ is always contained in the
same component, and $W_{r,2}$ in the other one.

Let us first check that for $j=1,2$,
$$
\Sigma_{r,j} \i U_{r,j} 
\leqno (3.34)
$$
where the $U_{r,j}$ are still the two connected components of
$\R^n \sm \Sigma_r$, defined below (3.21). By (3.20),
each point of $\Sigma_{r,j}$ is contained in some
$B(y,10^{-3}r)$, and then (3.34) follows from (3.23).
Let us deduce from this that
$$
W_{r,j} \i U_{r,j} \ \hbox{ for } j=1,2.
\leqno (3.35)
$$
Let $z\in W_{r,j}$ be given, pick any $\xi \in \Sigma_{r,j}$,
and let $w$ be the first point of $\Sigma_{r,j}$ when we run 
from $z$ to $\xi$ along $[z,\xi]$. Then $[z,w) \i W_{r,j}$
(because $W_{r,j}$ is a component of $\R^n \sm \Sigma_{r,j}$)
and hence $\dist([z,w),E) \geq 3C_0 \varepsilon$, by (3.32).
By (1.5), $\dist([z,w],E) \geq 2C_0 \varepsilon$, so
$[z,w]$ does not meet $\Sigma_r$, and $[z,w]$ is contained in
some $U_{r,i}$. But $w\in U_{r,j}$ by (3.34), 
so $i=j$ and (3.35) follows. 

Now we need to relate the descriptions that we get at 
different scales. Pick an origin $x_0 \in E$ and, for 
$0 < r  \leq r_0$, use (1.4) to choose $y(r) \in \Sigma_r$
such that $|y(r)-x_0| \leq C_0 \varepsilon r$.
Then denote by $P(r)$ and ${\cal G}(r)$ the plane and graph that
we get from (1.6)-(1.9), applied to $\Sigma_r$ and the point $y(r)$.
We claim that 
$$
d_{x_0,10^{-3}r}(P(r),P(s)) \leq C \varepsilon
\hbox{ when $0 < s \leq r \leq 9 s$.}
\leqno (3.36)
$$
Here $d$ is still the normalized Hausdorff distance defined in (1.2),
and the proof will be routine. First notice that 
${\cal G}(r)$ and $P(r)$ both go through $y(r)$ (by definitions).
If $a \in P(r) \cap B(x_0,10^{-3}r)$, we can use (1.9) to find
$z\in {\cal G}(r)$ such that $|z-a| \leq \Lambda_1 \varepsilon 
10^{-3}r$; this point lies in $\Sigma_r$ by (1.8),
and by (1.4) and (1.5) we can find $w \in \Sigma_s$
such that $|w-z| \leq 2 C_0 \varepsilon r$. Since
$w \in B(y(s),10^{-2} s)$, we can apply (1.8) again to see
that $w\in {\cal G}(s)$, and then its projection $b$
on $P(s)$ is such that $|b-w| \leq  \Lambda_1 \varepsilon 
10^{-3}s$. Thus we found $b\in P(s)$ such that
$|b-a| \leq C \varepsilon r$. The same argument shows that
for $b\in P(s) \cap B(x_0,10^{-3}r)$, we can find 
$a \in P(r)$ such that $|b-a| \leq C \varepsilon r$, and (3.36)
follows.

Denote by $e_n(r)$ the last vector of the basis that we chose
to describe $\Sigma_r$ near $y(r)$; thus $e_n(r)$ is orthogonal
to $P(r)$, and (3.36) shows that
$|e_n(r) \pm e_n(s)| \leq C \varepsilon$
when $0 < s \leq r \leq 9 s$; the sign $\pm$ depends on $r$ and $s$
through the orientations that we chose on $\Sigma_r$ and $\Sigma_s$.

Set $z_\pm(r) = x_0 + 10^{-4} r e_n(r)$ for $0 < r \leq r_0$; 
if follows (3.25) and (3.26) that $z_+(r) \in W_{r,1}$
and that $z_-(r) \in W_{r,2}$. But we can also look in the 
coordinates associated to $s$, use the fact that 
$|e_n(r) \pm e_n(s)| \leq C \varepsilon$, and get that
$z_\mp(r) \in W_{s,1}$ and $z_\pm(r) \in W_{s,2}$.
The signs are not important, but we want to know that
$z_+(r)$ and $z_-(r)$ lie in different $W_{s,j}$. 

We also need to know that when $r \geq 3s$,
each $W_{r,j}$ is contained in a $W_{s,i}$.
For each $z \in W_{r,j}$, (3.32) says that
$\dist(z,E) \geq 3C_0 \varepsilon r \geq 9C_0 s$.
Then by (3.30), $x\in W_{s,1} \cup W_{s,2}$.
That is, $W_{r,j} \i W_{s,1} \cup W_{s,2}$.
By (3.35), $W_{s,1}$ and $W_{s,2}$ are contained in different
connected components of $\R^n \sm \Sigma_s$, so
$W_{r,j}$, which is itself connected, is contained in a single 
$W_{s,i}$, as needed. Of course $i$ is unique, since 
$W_{s,1} \cap W_{s,2} = \emptyset$.

Let us now use all this to check that when 
$0 < 3s \leq r \leq r_0$,
$$
\hbox{$W_{r,1} \i W_{s,1}$ and $W_{r,2} \i W_{s,2}$, } 
\hbox{ or else $\, W_{r,1} \i W_{s,2}$ and $W_{r,2} \i W_{s,1}$.}
\leqno (3.37)
$$
Obviously it is enough to check this when $3s \leq r \leq 9 s$,
because after this we can use intermediate radii and iterate. 
Suppose first that $W_{r,1} \i W_{s,1}$. 
Then $z_+(r) \in W_{r,1} \cap W_{s,1}$,
and we observed earlier that then $z_-(r) \in W_{r,2} \cap W_{s,2}$.
In this case $W_{r,2} \i W_{s,2}$ (the other option, 
$W_{r,2} \i W_{s,1}$, is impossible), and we are happy.
If $W_{r,1}$ is not contained in $W_{s,1}$, then it is contained
in $W_{s,2}$, so $z_+(r) \in W_{r,1} \cap W_{s,2}$, then
$z_-(r) \in W_{r,2} \cap W_{s,1}$, and the only possibility is 
that $W_{r,2} \i W_{s,1}$; (3.37) follows.

We are now ready to show that for $0 < r \leq r_0$,
$W_{r,1}$ and $W_{r,2}$ lie in different components of
$\R^n \sm E$. First observe that by (3.32), 
$\dist(x,E) \geq 3C_0 \varepsilon r$ for $x\in W_{r,j}$,
so $W_{r,j}$ does not meet $E$ and, by connectedness, is contained
in a component of $\R^n \sm E$. If $W_{r,1}$ and $W_{r,2}$ lie in
the same component, we can find a (compact) curve in 
$\R^n \sm E$, that meets $W_{r,1}$ and $W_{r,2}$.
Since $\dist(\gamma,E) > 0$, (3.30) says that for $s$ small enough,
$\gamma \i W_{s,1} \cup W_{s,2}$. Since the $\gamma$ is connected
and the $W_{s,i}$ are contained in different components of 
$\R^n \sm \Sigma_s$ (by (3.35)), we see that $\gamma$ is contained
in a single $W_{s,i}$. But (3.37) precisely says that 
$W_{r,1}$ or $W_{r,2}$ is contained in the other $W_{s,j}$,
hence does not meet $W_{s,i}$. This contradiction proves (3.33).

We finally need to check that $\R^n \sm E$ has exactly two components.
Suppose instead that it has at least three components (we just 
excluded a single one). Let $z_1, z_2, z_3$ lie in different
components, and choose $r$ so small that for all $j$,
$\dist(z_j, E) > 6 C_0 \varepsilon r$; then (3.30) says that
all $z_j$ lie in $W_{r,1} \cup W_{r,2}$, at least two of them
lie in a same $W_{r,i}$, and this is a contradiction because
$W_{r,i}$ is connected and does not meet $E$.
This concludes our proof of Theorem 1.18.
\qed

\bigskip
REFERENCES
\smallskip 
\item {[AH]} P. Alexandroff and H. Hopf, Topologie.
Springer, Berlin, 1935.
\smallskip 
\item {[DS]}
Guy David and  Stephen Semmes,
\underbar{Analysis of and on uniformly rectifiable sets},
Mathematical Surveys and Monographs 38, American Mathematical 
Society, 
Providence, RI, 1993. xii+356 pp.
\smallskip
\item {[DT]} G. David and T. Toro,
Reifenberg parameterizations for sets with holes.
Memoirs of the American Mathematical Society,
Vol. 215, 2012, Number 1012.
\smallskip
\item {[F]} H. Federer, \underbar{Geometric measure theory}, 
Grundlehren der Mathematishen Wissenschaf-ten 
153, Springer Verlag 1969.
\smallskip
\item {[R]} E.\ R.\ Reifenberg, Solution of the Plateau problem
for $m$-dimensional surfaces of varying topological type, 
Acta Mathematica {\bf 104} (1960), 1--92.
\smallskip
\item {[S]} H. Samelson,
Orientability of hypersurfaces in $\R^{n}$. 
Proc. Amer. Math. Soc. 22, 1969, 301Ð302. 
\smallskip
\item {[St]} E. M. Stein, \underbar{Singular integrals and 
differentiability properties of functions},
Princeton university press 1970.
\smallskip

\bsi\bsi
Guy David, 
\pari
Univ Paris-Sud, 
\pari
Laboratoire de math\'{e}matiques UMR-8628, 
\pari
Orsay F-91405, France
\pari
and 
\pari
Institut Universitaire de France

\bye